\documentclass[11pt]{article} 
\title{{\bf Pressure Rigidity of Three Dimensional 
Contact Anosov Flows}} 
\author{Yong Fang} 
\date{U.M.R. 
7501 du C.N.R.S, 
Institut de Recherche Math\'ematique Avanc\'ee, 7 rue 
Ren\'e Descartes, 67084 Strasbourg C\'edex, France\\ (e-mail : fang@math.u-strasbg.fr)} 

\usepackage{amsmath,amsthm,amsfonts,amssymb,amscd,amstext} 
\chardef\bslash=`\\

\theoremstyle{definition}

\theoremstyle{remark}

\begin{document} 
\maketitle 
\renewcommand{\sectionmark}[1]{} 
{\bf Abstract} --- {\it Let $\phi_t$ be a three dimensional contact Anosov flow. Then we prove 
that its cohomological pressure coincides with its metric entropy if and only if $\phi_t$ is $C^\infty$ flow 
equivalent to 
a special time change of a three dimensional algebraic contact Anosov flow.}\\\\
{\bf 1. Introduction}\\\\
Let $M$ be a $C^\infty$-closed manifold. A $C^\infty$-flow, $\phi_t$, generated by a non-singular vector field $X$ 
on $M$ is called 
{\it Anosov}, if there exists a $\phi_t$-invariant continuous splitting of the tangent bundle
$$TM =\mathbb{R}X\oplus E^{+}\oplus E^{-},$$
a Riemannian metric on $M$ and two positive numbers $a$ and $b$, such that 
$$\parallel D\phi_{-t}
(u^{+})\parallel\leq a\cdot e^{-bt}\parallel u^{+}\parallel, \  \forall \  u^{+}\in E^{+},\  \forall \  t > 0,$$
and
$$\parallel D\phi_{t}
(u^{-})\parallel\leq a\cdot e^{-bt}\parallel u^{-}\parallel, \  \forall \  u^{-}\in E^{-},\  \forall \  t > 0.$$
The continuous distributions 
$E^{+}$ and $E^{-}$ are called 
respectively the strong unstable and strong stable distributions of $\phi_t$. They are both 
integrable to continuous foliations with $C^\infty$ leaves (see {\bf [HK]}).
 
The {\it canonical $1$-form} of $\phi_t$ is by definition the continuous $1$-form 
on $M$ such that $\lambda(X)=1$ and $\lambda(E^\pm)=0$. It is easily seen that $\lambda$ is $\phi_t$-invariant. By 
definition, $\phi_t$ is said to be a {\it contact Anosov} flow 
if $\lambda$ is $C^\infty$ and there exists $n\in\mathbb{N}$ such that $\lambda\wedge (\wedge^n d\lambda)$ is a 
volume form on $M$. It is easy to see that contact Anosov flows are contact 
in the classical sense (see {\bf [Pa]}).    

Let $\phi_t$ be a contact Anosov flow 
on a closed manifold $M$ of dimension $2n+1$. Then $\lambda\wedge(\wedge^n d\lambda)$ is a 
$\phi_t$-invariant volume form. Denote by $\nu$ 
the $\phi_t$-invariant probability measure determined by this volume form. Then the 
measure-theoretic 
entropy of $\phi_t$ with respect to $\nu$ is said to be the metric entropy of $\phi_t$ and is denoted by 
$h_\nu(\phi_t)$. 

Denote by $\mathcal{M}(\phi_t)$ the set of $\phi_t$-invariant probability measures. For each 
$C^\infty$ function $f$ on $M$, the {\it topological pressure} of 
$\phi_t$ with respect to $f$ is dy definition the following number
$$P(\phi_t, f) = sup_{\mu\in\mathcal{M}(\phi_t)}\{h_\mu(\phi_t)+\int_M f d\mu\}, $$
where $h_\mu(\phi_t)$ denotes the metric entropy of $\phi_t$ with respect to $\mu$.  
It is well known (see {\bf [HK]} and 
{\bf [BR]}) that there exists a unique $\phi_t$-invariant probability measure $\mu$ such that 
$$ P(\phi_t, f) = h_{\mu}(\phi_t)+\int_M f d\mu.$$
This measure $\mu$ is said to be the {\it Gibbs measure} of $\phi_t$ with respect to $f$. For 
example, the classical Bowen-Margulis of $\phi_t$ is just the Gibbs measure of 
$\phi_t$ with respect to the zero function. Recall also that $P(\phi_t, 0)$ is just the topological entropy of 
$\phi_t$ denoted by $h_{top}(\phi_t)$.  

Let $f'$ be another 
$C^\infty$ function on $M$ with Gibbs measure $\mu'$. Then we can 
prove (see {\bf [HK]}) that $\mu=\mu'$ if and only if there exist a $C^\infty$ 
funcion $H$ and a constant $c$ such that
$$ f= f'+X(H)+c.$$
Denote by $H^1(M, \mathbb R)$ the first cohomology group of $M$. 
For $[\beta]\in H^1(M, \mathbb{R})$ we denote by $P(\phi_t, [\beta])$ the topological pressure of $\phi_t$ 
with respect to $\beta(X)$. 
For each smooth function $g$ and $\mu\in \mathcal{M}(\phi_t)$ we have
$$\int_{M}X(g) d\mu=0.$$
So $P(\phi_t, [\beta])$ is independent of the closed $1$-form chosen in the cohomological class $[\beta]$. 
We call the Gibbs measure of $\phi_t$ with 
respect to $\beta(X)$ that of $\phi_t$ with respect to $[\beta]$. We define 
$$ P(\phi_t)= inf_{[\beta]\in H^1(M, \mathbb R)}\{ P(\phi_t, [\beta])\},$$
which is said to be the {\it cohomological pressure} of $\phi_t$. 
This notion was firstly defined by R. Sharp in {\bf [Sh]}. By the equivalence of $(ii)'$ and $(iii)$ of Theorem one 
in {\bf [Sh]}, we know that there exists a unique element $[\alpha]$ in $H^1(M, \mathbb R)$ such that 
$$ P(\phi_t, [\alpha])= P(\phi_t).$$
We call this cohomology class $[\alpha]$ the {\it Gibbs class} of $\phi_t$. 

For each element $[\beta]$ in $H^1(M, \mathbb R)$, it is easy to see that  
$$ \int_{M}\beta(X)\lambda\wedge(\wedge^n d\lambda)= \int_M\beta\wedge(\wedge^n d\lambda)=0.$$
So $ \int_{M}\beta(X) d\nu=0.$ Thus we have
$$ h_{top}(\phi_t)\geq P(\phi_t)\geq h_\nu(\phi_t).$$
In general these inequalities are strict.

Let $N$ be a $C^\infty$ closed negatively curved manifold. Denote by $\phi_t$ its geodesic flow and 
by $\mu$ the Bowen-Margulis measure of $\phi_t$. Since 
$\mu$ is invariant under the flip map, then for each $[\alpha]$ in $H^1(SN, \mathbb R)$ we have
$$\int_{SN}\alpha(X)d\mu=0,$$
where $SN$ denotes the unitary bundle of $N$ (see {\bf [Pa]}). So we get
$$P(\phi_t)\geq inf_{[\alpha]\in H^1(SN,\mathbb R)}\{h_\mu(\phi_t)+\int_{SN}\alpha(X)d\mu\}= h_\mu
(\phi_t)= h_{top}(\phi_t).$$
We deduce that $P(\phi_t)= h_{top}(\phi_t)$. So for the geodesic flows of negatively curved 
manifolds the cohomological pressure coincides with the topological entropy.\\\\
{\bf 2. Rigidity in the case of dimension three}\\\\    
The classical examples of three dimensional contact Anosov flows are constructed as following. 
Denote by $\Gamma$ 
a uniform lattice in $\widetilde{SL(2, \mathbb{R})}$. Then on the quotient manifold 
$\Gamma\diagdown\widetilde{SL(2, \mathbb{R})}$, the following flow is 
contact Anosov and is said to be {\it algebraic}. 
$$\phi_t: \Gamma\diagdown\widetilde{SL(2, \mathbb{R})}\to \Gamma\diagdown\widetilde{SL(2, \mathbb{R})},$$
$$ \Gamma\cdot A\to \Gamma\cdot (A\cdot {\rm exp}(t \left (\begin{array}{cccc}
1 & 0 \\
0 & -1
\end{array}
\right))).$$
Up to a constant change 
of time scale and finite covers, such a flow is just the geodesic flow of a certain closed hyperbolic surface.

Quite recently, P. Foulon constructed in {\bf [Fo]} plenty 
of surgerical three dimensional 
contact Anosov flows. These examples make the study of three dimensional 
contact Anosov flows very interesting.

Now suppose that $\phi_t$ is a three dimensional contact Anosov flow on a closed manifold $M$. A 
{\it special time change} of 
$\phi_t$ is by definition the flow of $\frac{X}{a-\alpha(X)}$, where $a>0$ and 
$\alpha$ denotes a $C^\infty$ closed $1$-form on $M$ such that $a-\alpha(X)>0.$\\\\ 
{\bf Lemma 2.1.} {\it Let $\phi_t$ be a 
three dimensional contact Anosov flow and $\psi_t$ be a smooth time change of $\phi_t$. Then 
$\psi_t$ is contact iff it is a special time change of $\phi_t$.}\\\\
{\bf Proof.} Suppose that $\psi_t$ is contact. Then its canonical $1$-form $\bar\lambda$ is $C^\infty$. 
Denote by $Y$ the generator of $\psi_t$ and suppose that $Y= fX$. Then we have 
$$f\cdot\bar\lambda(X) = \bar \lambda(Y)=1.$$
So we have $Y= \frac{X}{\bar\lambda(X)}$. Since $\bar\lambda$ is $\psi_t$-invariant, then 
$$ i_Y d\bar\lambda=0.$$
We deduce that
$$ \mathcal{L}_X d\bar\lambda= di_Xd\bar\lambda+ i_X dd\bar\lambda=0.$$
So $d\bar\lambda$ is $\phi_t$-invariant. Denote by $\lambda$ the canonical $1$-form of $\phi_t$. Then by the 
ergodicity of $\phi_t$ (see {\bf [An]}) there 
exists a constant $b$ such that
$$ \lambda\wedge d\bar\lambda= b\cdot \lambda\wedge d\lambda.$$
So there exists a $C^\infty$ closed $1$-form $\alpha$ such that
$$ \bar \lambda= b\cdot\lambda +\alpha.$$
If $b$ is non-positive, then $\alpha(X)>0$ on $M$. Denote by $Z$ the field $\frac{X}{\alpha(Z)}$. Then $\alpha$ is 
easily seen to be the canonical $1$-form of $\phi^Z_t$. So by 
{\bf [Pl]}, $\phi_t$ admits a global section, which is 
absurd for a contact flow. We deduce that $b>0$. Thus $\psi_t$ is a special time change of $\phi_t$. 

Suppose that $\psi_t$ is a special time change of $\phi_t$, i.e. $Y= \frac{X}{a-\alpha(X)}.$ 
It is easily seen that $a\lambda-\alpha$ is the canonical $1$-form of $Y$. We have 
$$(a\lambda-\alpha)\wedge d(a\lambda-\alpha)= a^2\cdot \lambda\wedge d\lambda- a\cdot \alpha\wedge d\lambda.$$
By integrating this form, we see that $(a\lambda-\alpha)\wedge d(a\lambda-\alpha)\not\equiv 0.$ 
Then we deduce by {\bf [HuK]} that 
this three form is nonwhere zero, i.e. $\psi_t$ is contact.  $\square$

$\ $

In the quite elegant paper {\bf [Ka]}, A. Katok proved the following\\\\ 
{\bf Theorem 2.1.} (A. Katok) {Let $\Sigma$ be a $C^\infty$ closed surface of negative curvature. Then its 
topological entropy coincides with its metric entropy, if and only if $\Sigma$ is of constant negative curvature.}

$\ $

Then in {\bf [Fo1]}, the following generalization was established by using geometric constructions.\\\\
{\bf Theorem 2.2.} (P. Foulon) {\it Let $\phi_t$ be a three dimensional contact Anosov flow. Then its 
topological entropy coincides with its metric entropy if and only if 
it is, up to a constant change of time scale, $C^\infty$ flow equivalent to a three dimensional 
algebraic contact Anosov flow.}

$\ $

Now we generalize this Theorem to the case of cohomological pressure. Let us prove firstly the following\\\\
{\bf Lemma 2.2.} {\it Let $\phi_t$ be a three dimensional algebraic contact Anosov flow. Then for each 
special time change $\psi_t$ of $\phi_t$, there exists an element $[\alpha]$ in $H^1(M, \mathbb{R})$ such that the 
Gibbs measure of $\psi_t$ 
with respect to $[\alpha]$ is Lebesgue. In particular the cohomological 
pressure of $\psi_t$ coincides with its 
metric entropy.}\\\\
{\bf Proof.} Denote by $h$ the topological entropy of $\phi_t$. 
Since the Anosov splitting of $\phi_t$ is $C^\infty$, then we can 
find a $C^\infty$ nonwhere vanishing section $Y^+$ of 
$E^+$ and define for any $x\in M$ and any $t\in\mathbb{R}$,
$$\lambda_t(x) = \frac{(\phi_t)_\ast Y^+_x}{Y^+_{\phi_t(x)}}.$$
Define also
$$ \phi^+ = -\frac{\partial}{\partial t}\mid_{t=0}\lambda_t(\cdot).$$
Then it is easy to see 
that the Gibbs measure of $\phi_t$ with respect to $\phi^+$ is $\nu$ (see {\bf [BR]}). In addition, we have 
$$P(\phi_t, \phi^+)=0.$$
Since the Bowen-Margulis measure of $\phi_t$ is Lebesgue, then there exists a smooth function $f$ and 
a contant $c$ such that
$$ \phi^+ = X(f) +c.$$
So we have 
$$0= P(\phi_t, \phi^+)= P(\phi_t, X(f)+c) = h+c.$$
We deduce that $\phi^+= X(f)-h.$ 

Suppose that $\psi_t$ is generated by the field $\bar X= \frac{X}{a-\alpha(X)}$. Then by {\bf [LMM]}, it is easy to 
see that  
$$\bar E^+= \{ u^+ +\frac{\alpha(u^+)}{a-\alpha(X)}\cdot X\mid u^+\in E^+\}.$$
Suppose that $\psi_t(\cdot)= \phi_{\beta(t, \cdot)}(\cdot)$ and define $\bar Y^+= Y^+
+\frac{\alpha(Y^+)}{a-\alpha(X)}X.$ Then
$$\bar\lambda_t(x)= \frac{(\psi_t)_\ast\bar Y^+_x}{\bar Y^+_{\psi_t(x)}} = \lambda_{\beta(t, x)}(x).$$
Since $\dot\beta(0, \cdot)= \frac{1}{a-\alpha(X)}$, then 
$$\bar \phi^+= \frac{\phi^+}{a-\alpha(X)}.$$
So we get 
$$\bar\phi^+= \bar X(f) - \frac{h}{a-\alpha(X)}.$$
In addition, we observe easily that
$$(-\frac{h}{a}\cdot\alpha)(\bar X) = \frac{h}{a}- \frac{h}{a-\alpha(X)}.$$
So the Gibbs measure of $\psi_t$ with respect to $[-\frac{h}{a}\cdot\alpha]$ is Lebesgue.  $\square$

$\ $

Now we establish the following extension of Theorem $2.2.$\\\\ 
{\bf Theorem 2.3.} {\it Let $\phi_t$ be a three dimensional contact Anosov flow defined on a closed 
manifold $M$. Then its cohomological pressure coincides with its metric entropy if and only if 
$\phi_t$ is $C^\infty$ flow 
equivalent to 
a special time change of a three dimensional algebraic contact Anosov flow.}\\\\
{\bf Proof.} Denote by $[\alpha]$ the Gibbs class of $\phi_t$. Since $\int_M\alpha(X)d\nu=0$ and by assumption 
$$ h_\nu(\phi_t)= P(\phi_t)= P(\phi_t, [\alpha]),$$
then 
$$h_\nu(\phi_t)+\int_M\alpha(X)d\nu= P(\phi_t, [\alpha]),$$
i.e. $\nu$ is the Gibbs measure of $\phi_t$ with respect to $[\alpha]$. Denote $h_\nu(\phi_t)$ 
by $h$. Then by the variational principle, we have for $\forall\  \mu\in \mathcal M(\phi_t)$ and $\mu\not= \nu$,
$$ h>h_\mu(\phi_t)+\int_M\alpha(X) d\mu\geq \int_M\alpha(X) d\mu.$$
Then by {\bf [Gh1]}, there exists a smooth function $f$ such that  
$$h>(\alpha+ df)(X).$$
Without loss of generality, we replace $\alpha$ by $\alpha+df$. 

Up to finite covers, we suppose that $E^+$ and $E^-$ are both orientable. In {\bf [HuK]}, 
it is proved that the strong stable and instable 
distributions of a three dimensional contact Anosov flow are $C^{1, Zyg}$. So we can find a $C^{1, Zyg}$ 
nonwhere vanishing 
section $Y^+$ of $E^+$ and define the functions $\lambda_t(\cdot)$ and $\phi^+$ as in the proof of 
Lemma $2.3$. Then the Gibbs measure of $\phi_t$ 
with respect to $\phi^+$ is lebesgue. 

Define $\bar X=\frac{X}{h-\alpha(X)}$ and denote by $\psi_t$ its flow. Then we define as before $\bar Y^+$ and 
$\bar\phi^+$. So by the same arguments, 
$$ \bar\phi^+= \frac{\phi^+}{h-\alpha(X)}.$$
Since the Gibbs measure of $\phi_t$ 
with respect to $\alpha(X)$ is Lebesgue, then there exists a smooth function 
$g$ and a constant $c$ such that
$$ \phi^+= \alpha(X)+ X(g)+c.$$
So we have 
$$ 0= P(\phi_t,\phi^+)= P(\phi_t, [\alpha])+c= h+c.$$
We deduce that 
$$ \bar\phi^+= -1+ \bar X(g),$$
i.e. the Bowen-Margulis measure of $\psi_t$ is Lebesgue. Then by Theorem $2.2$, $\psi_t$ is $C^\infty$ 
flow equivalent to a three dimensional algebraic contact Anosov flow. We deduce from 
Lemma $2.2$ that $\phi_t$ is $C^\infty$ flow equivalent to a special time change of such a flow.  $\square$

$\ $

As mentioned above, we know by {\bf [HuK]} that the 
Anosov splitting of a three dimensional contact Anosov flow is always $C^{1, Zyg}$. 
In addition by Theorem $4.6$ of {\bf [Gh2]}, we know that up to finite covers, a three 
dimensional contact Anosov flow with $C^{1, lip}$ splitting is $C^\infty$ flow equivalent to 
a special time change of the geodesic flow of a closed hyperbolic surface (see also 
{\bf [Gh]}, {\bf [BFL]} and {\bf [HuK]}). Thus by combining these classical results with our previous result, we obtain 
the following \\\\
{\bf Theorem 2.4.} {\it Let $\phi$ be a three dimensional contact Anosov flow. Then 
its Anosov splitting is $C^{1, Zyg}$ and the following conditions are 
equivalent : \\
$(1)$ The cohomological pressure of $\phi_t$ is equal to its metric entropy.\\
$(2)$ The Anosov splitting of $\phi_t$ is $C^{1, lip}$.\\
$(3)$ Up to a constant 
change of time scale and finite covers, $\phi_t$ is $C^\infty$ flow equivalent to 
a special time change of the geodesic 
flow of a closed 
hyperbolic surface.}\\\\          
{\bf Acknowledgements.} The author would like to thank F. Ledrappier and P. Foulon 
for interesting discussions.\\\\

{\bf References}

$\ $

{\small

{\bf [An]} V. D. Anosov, Geodesic flows on closed Riemannian manifolds with negative curvature, 
{\it Proc. Inst. Steklov} 90 (1967) 1-235.

{\bf [BFL]} Y. Benoist, P. Foulon and F. Labourie, Flots d'Anosov \`a distributions 
stable et instable diff\'erentiables, {\it J. Amer. Math. Soc.} 5 (1992) 33-74.

{\bf [BR]} R. Bowen and D. Ruelle, The ergodic theory of Axiom A flows, {\it Invent. Math.} 29 (1975) 181-202.

{\bf [Fo]} P. Foulon, {\it personal communications.}

{\bf [Fo1]} P. Foulon, Entropy rigidity of Anosov flows in dimension three, {\it Ergod. Th. and Dynam. Sys.} 21 
(2001) 1101-1112.

{\bf [Gh]} \'E. Ghys, Flots d'Anosov dont les feuilletages stables sont 
diff\'erentiables, {\it Ann. Scient. \'Ec. Norm. Sup. (4)} 20 (1987) 251-270.

{\bf [Gh1]} \'E. Ghys, Codimension one Anosov flows and suspensions, {\it Lecture Notes in 
Mathematics} 1331 (1988) 59-72.

{\bf [Gh2]} \'E. Ghys, Rigidit\'e diff\'erentiable des groupes fuchsiens, 
{\it Pub. IH\'ES } 78, (1993) 163-185.
 
{\bf [HK]} B. Hasselblatt and A. Katok, Introduction to the modern theory of dynamical 
systems, {\it Encyclopedia of Mathematics and its Applications, vol 54.} 1995. 

{\bf [HuK]} S. Hurder and A. Katok, Differentiability, rigidity and 
Godbillon-Vey classes for Anosov flows, {\it Pub. I.H.\'E.S.} 72 
(1990) 5-61.

{\bf [Ka]} A. Katok, Entropy and closed geodesics, {\it Ergod. Th. and Dynam. Sys.} 2 (1982) 339-367.

{\bf [LMM]} R. de la llave, J. Marco and R. Moriyon, Canonical perturbation theory of 
Anosov systems and regularity results for Livsic cohomology equation, {\it Ann. Math} 
123(3) (1986) 537-612.

{\bf [Pa]} G. P. Paternain, Geodesic flows, {\it Progress in Mathematics.}

{\bf [Pl]} J. F. Plante, Anosov flows, {\it Amer. J. Math.} 94 (1972) 729-754.

{\bf [Sh]} R. Sharp, Closed orbits in homology classes for Anosov flows, {\it Ergod. Th. and Dynam. Sys.} 13 (1993) 
387-408.}

\end{document}